\pdfoutput=1
\RequirePackage{ifpdf}
\ifpdf 
\documentclass[pdftex]{sigma}
\else
\documentclass{sigma}
\fi

\numberwithin{equation}{section}

\newtheorem{Theorem}{Theorem}[section]
\newtheorem{Conjecture}[Theorem]{Conjecture}

\newcommand{\di}{\mathop{}\!\mathrm{d}}
\newcommand{\dist}{\mathsf{d}}

\newcommand{\meas}{\mathfrak{m}}
\DeclareMathOperator{\RCD}{RCD}

\begin{document}
\allowdisplaybreaks

\newcommand{\arXivNumber}{2002.08612}

\renewcommand{\thefootnote}{}

\renewcommand{\PaperNumber}{021}

\FirstPageHeading

\ShortArticleName{Collapsed Ricci Limit Spaces as Non-Collapsed ${\rm RCD}$ Spaces}

\ArticleName{Collapsed Ricci Limit Spaces\\ as Non-Collapsed $\boldsymbol{{\rm RCD}}$ Spaces\footnote{This paper is a~contribution to the Special Issue on Scalar and Ricci Curvature in honor of Misha Gromov on his 75th Birthday. The full collection is available at \href{https://www.emis.de/journals/SIGMA/Gromov.html}{https://www.emis.de/journals/SIGMA/Gromov.html}}}

\Author{Shouhei HONDA}

\AuthorNameForHeading{S.~Honda}

\Address{Mathematical Institute, Tohoku University, Sendai 980-8578, Japan}
\Email{\href{shouhei.honda.e4@tohoku.ac.jp}{shouhei.honda.e4@tohoku.ac.jp}}

\ArticleDates{Received February 20, 2020, in final form March 25, 2020; Published online April 01, 2020}

\Abstract{In this short note we provide several conjectures on the regularity of measured Gromov--Hausdorff limit spaces of Riemannian manifolds with Ricci curvature bounded below, from the point of view of the synthetic treatment of lower bounds on Ricci curvature for metric measure spaces.}


\Classification{53C20; 53C21; 53C23}

\begin{flushright}
\begin{minipage}{65mm}
\it Dedicated to Professor Mikhail Gromov\\ on the occasion of his 75th birthday.
\end{minipage}
\end{flushright}

\renewcommand{\thefootnote}{\arabic{footnote}}
\setcounter{footnote}{0}

\section{Convergence theory of Riemannian manifolds}

Let $\big(M^n_i, g_i\big)$ be a sequence of $n$-dimensional complete Riemannian manifolds with the pointed Gromov--Hausdorff (pGH) convergence:
\begin{gather}\label{111}
(M^n_i, \dist_{g_i}, x_i) \stackrel{\mathrm{pGH}}{\longrightarrow} (X, \dist, x)
\end{gather}
to some pointed proper metric space $(X, \dist, x)$, where a metric space is said to be proper if any bounded closed subset is compact.
Then the convergence theory of Riemannian manifolds states that under suitable curvature restriction:
\begin{enumerate}\itemsep=0pt
\item[$(\star)$] Establish regularity results on $(X, \dist)$.
\item[$(\star \star)$] Find relationships between $M^n_i$ and $X$.
\end{enumerate}

It is worth pointing out that Gromov's cerebrated precompactness theorem states that a~sequence of pointed proper metric spaces $(X_i, \dist_i, x_i)$ has a pGH-convergent subsequence if and only if open balls $B_R(x_i)$ are uniformly metric doubling for any fixed $R>0$, that is, for any~$\epsilon>0$ there exists $k \in \mathbb{N}$ such that for any $i \in \mathbb{N}$ there exists a collection of points $y_{i, 1}, \ldots, y_{i, k} \in B_R(x_i)$ such that $B_R(x_i) \subset \bigcup_{j=1}^kB_{\epsilon}(y_{i, j})$ holds.

In this short note we provide several conjectures related to the problem~$(\star)$ under lower Ricci curvature bounds. Before introducing them, let us start to discuss the case of sectional curvature shortly in order to clarify the difference from the case of Ricci curvature.

Note that metric (metric measure, respectively) spaces with a lower bound of sectional (Ricci, respectively) curvature we will discuss below satisfy this uniform metric doubling property because of the Bishop--Gromov inequality. Thus such a sequence always has a pGH-convergent subsequence, which shows us that the setting~(\ref{111}) appears naturally in the following settings.

\section{Lower bound on sectional curvature}
Let us consider (\ref{111}) in the case when the sectional curvature of $\big(M^n_i, g_i\big)$ is bounded below by a constant $K \in \mathbb{R}$:
\begin{gather}\label{2}
\mathrm{Sec}_{M^n_i}^{g_i} \ge K.
\end{gather}
Then $(X, \dist)$ is a $k$-dimensional Alexandrov space of curvature bounded below by $K$ for some $k \in [0, n] \cap \mathbb{N}$, which is a direct consequence of the stability of Alexandrov spaces with respect to the pGH convergence proved in~\cite{BGP}. In particular nice geometric properties of~$(X, \dist)$ are carried from Alexandrov geometry which gives the best framework on the synthetic treatment of lower bounds on sectional curvature for metric spaces. For example any point~$p$ of~$X$ has a neighbourhood $U_p$ of $p$ which is homeomorphic to the tangent cone at~$p$, which is proved in \cite{BGP, P} (see also \cite{K}). This gives a geometric answer to the problem~$(\star)$.
The fibration theorem proved in \cite{Y} also gives a geometric answer to the problem $(\star \star)$.

\section{Lower bound on Ricci curvature}
Next let us consider (\ref{111}) in the case:
\begin{gather}\label{3}
\mathrm{Ric}_{M^n_i}^{g_i} \ge K(n-1).
\end{gather}
It is trivial that (\ref{2}) implies (\ref{3}).
After passing to a subsequence by \cite{CheegerColding1, F}, with no loss of generality we can assume that the pointed measured Gromov--Hausdorff (pmGH) convergence holds:
\begin{gather}\label{1}
\left(M^n_i, \dist_{g_i}, x_i, \frac{\mathrm{vol}_{g_i}}{\mathrm{vol}_{g_i}B_1(x_i)} \right) \stackrel{\mathrm{pmGH}}{\longrightarrow} (X, \dist, x, \meas)
\end{gather}
for some Borel measure $\meas$ on $X$.
Then $(X, \dist, \meas)$ is so-called a \textit{Ricci limit space}.

The structure theory of Ricci limit spaces is established in \cite{CheegerColding1, CheegerColding2, CheegerColding3}. For example, the Laplacian $\Delta$ on $(X, \dist, \meas)$ is well-defined via the rectifiablity as a metric measure space. Moreover it is proved in \cite{CheegerColding3} that the spectrums behave continuously with respect to the convergence (\ref{1}) if $(X, \dist)$ is compact, which confirms a conjecture raised in \cite{F}. This gives an analytic answer to the problem $(\star \star)$.

On the geometric side, in general, a similar fibration result as in \cite{Y} is not satisfied in this setting. A counterexample can be found in \cite{An}. However we know that if $(X, \dist)$ is compact with no singular set, and (\ref{1}) is non-collapsed (as explained below), then $X$ is homeomorphic to $M^n_i$ for any sufficiently large $i$, which is proved in \cite{CheegerColding1}. This gives a geometric answer to the problem $(\star \star)$.

Let us consider the problem $(\star)$. It is proved in \cite{CheegerColding1} that the same geometric property as in the previous section also holds for regular points if the sequence (\ref{1}) is \textit{non-collapsed} whose definition is to satisfy $\meas=a \mathcal{H}^n$ for some $a \in (0, \infty)$, where $\mathcal{H}^n$ is the $n$-dimensional Hausdorff measure. That is, if $\meas=a \mathcal{H}^n$, then any $n$-dimensional regular point~$p$ of~$X$ has a~neighbourhood~$U_p$ of $p$ which is homeomorphic to $\mathbb{R}^n$.\footnote{We say that a~point~$p$ is $n$-dimensional regular if any tangent cone at~$p$ is isometric to $\big(\mathbb{R}^n, \dist_{\mathbb{R}^n}, 0_n\big)$. In general tangent cones at a point are not unique even in the non-collapsed setting. More strongly, there exists a $5$-dimensional non-collapsed Ricci limit space with a~base point~$p$ such that there exist two tangent cones at~$p$ which are not homeomorphic to each other. See~\cite{CN1}. Thus it is hard to find nice topological results around singular points, which is very different from the case of sectional curvature~(\ref{2}). See also~\cite{M1}.} This gives a geometric answer to the problem~$(\star)$. See also~\cite{CJN} for a recent development along this direction.

However if the sequence (\ref{1}) is collapsed, then such a nice geometric property is unknown. Although one of the central topics in this story is to develop the local structure theory in general situation, it is still very hard.

In connection with these observations,
it is also interesting to ask:
\begin{enumerate}\itemsep=0pt
\item[$(\spadesuit)$] When can we find other non-collapsed sequence of $k$-dimensional complete Riemannian manifolds $\big(\hat{M}^k_i, \hat{g}_i\big)$ with Ricci curvature bounded below by a constant such that
\begin{gather}\label{15}
\left(\hat{M}^k_i, \dist_{\hat{g}_i}, \hat{x}_i, \frac{\mathrm{vol}_{\hat{g}_i}}{\mathrm{vol}_{\hat{g}_i}B_1(\hat{x}_i)} \right) \stackrel{\mathrm{pmGH}}{\longrightarrow} (X, \dist, x, \meas)
\end{gather}
holds as a non-collapsed sequence (even if the sequence (\ref{1}) is collapsed)?
\end{enumerate}
In general this question $(\spadesuit)$ has a negative answer. A counterexample can be found in \cite{CheegerColding1} as a~``metric horn'' which will be discussed later. See also \cite{H, M2, M3} for other examples along this direction. Let us emphasize that if (\ref{15}) holds (as a non-collapsed sequence), then $\meas=b \mathcal{H}^k$ holds for some $b \in (0, \infty)$.

We are now in a position to provide the first conjecture:
\begin{Conjecture}\label{4}
If $\meas=b\mathcal{H}^k$ holds for some $b, k \in (0, \infty)$, then $\big(X, \dist, \mathcal{H}^k\big)$ is a non-collapsed $\RCD(K(n-1), k)$ space.
\end{Conjecture}
The synthetic condition of lower bounds on Ricci curvature for metric measure spaces, $\RCD(K, N)$ condition, is explained in the next section.

Let us compare Conjecture~\ref{4} with the following conjecture raised in~\cite{CheegerColding1}:
\begin{Conjecture}[Cheeger--Colding]\label{200}
It holds that there exists $m \in [k, n] \cap \mathbb{N}$ such that for all $y \in X$,
\begin{gather*}
\frac{\mathcal{H}^k(B_r(y))}{\mathrm{Vol}_{m}(r)} \downarrow \qquad (r \uparrow \infty)
\end{gather*}
holds, where $k=\dim_{\mathcal{H}}(X, \dist)$ is the Hausdorff dimension of $(X, \dist)$ and $\mathrm{Vol}_m(r)$ denotes the volume of a ball of radius $r$ in the $m$-dimensional space form whose sectional curvature is equal to $-1$.
\end{Conjecture}
In connection with Conjecture \ref{200}, it is natural to ask:
\begin{enumerate}\itemsep=0pt
\item[$(\heartsuit)$] Is $\big(X, \dist, \mathcal{H}^{k}\big)$ an $\RCD(K_1, N_1)$ space for some $K_1 \in \mathbb{R}$ and some $N_1 \in [1, \infty]$? Here $k=\dim_{\mathcal{H}}(X, \dist)$.
\end{enumerate}
Because if this question $(\heartsuit)$ has a positive answer for some suitable $K_1$, $N_1$, then Conjecture~\ref{200} holds by the Bishop--Gromov inequality in the~$\RCD$ theory (note that then it is necessary to satisfy $N_1 \ge k$). However this question $(\heartsuit)$ has a negative answer for a metric horn. See the next section.

\section{Synthetic treatment of lower bound on Ricci curvature}
A triple $(X, \dist, \meas)$ is said to be a metric measure space if $(X, \dist)$ is a complete separable metric space and $\meas$ is a Borel measure on $X$ with full support. The pioneer papers \cite{LottVillani, Sturm06a, Sturm06b} define $\mathrm{CD}(K, N)$ conditions for metric measure spaces from the point of view of optimal transportation theory (see also~\cite{Villani}). Roughly speaking, a metric measure space $(X, \dist, \meas)$ is said to be a~$\mathrm{CD}(K, N)$ space if the Ricci curvature is bounded below by $K \in \mathbb{R}$, and the dimension is bounded above by $N \in [1, \infty]$.
After that, adding Riemannian structure to $\mathrm{CD}(K, N)$ spaces, the definition of $\RCD(K, N)$ spaces is introduced in \cite{AmbrosioGigliSavare14, ErbarKuwadaSturm, Gigli13}. Let us give an equivalent definition of~$\RCD(K, N)$ spaces under assuming a bit of knowledges on the Sobolev space $H^{1, 2}(X, \dist, \meas)$.

A metic measure space $(X, \dist, \meas)$ is an $\RCD(K, N)$ \textit{space} for some $K \in \mathbb{R}$ and some $N \in [1, \infty]$ if the following four conditions hold:
\begin{itemize}\itemsep=0pt
\item{(Volume growth condition)} There exist $C>1$ and $x \in X$ such that $\meas (B_r(y)) \le Ce^{C\dist (x, y)^2}$ holds for any $y \in X$ and any $r>0$.
\item{(Riemannian structure)} The Sobolev space $H^{1, 2}=H^{1, 2}(X, \dist, \meas)$ is a Hilbert space. In particular for all $f_i \in H^{1, 2}$ $(i=1, 2)$,
\begin{gather*}
\langle \nabla f_1, \nabla f_2\rangle:=\lim_{t \to 0}\frac{|\nabla (f_1+tf_2)|^2-|\nabla f_1|^2}{2t} \in L^1(X, \meas)
\end{gather*}
is well-defined, where $|\nabla f_i|$ denotes the minimal relaxed slope of $f_i$.
\item{(Sobolev-to-Lipschitz property)} Any function $f \in H^{1, 2}$ satisfying $|\nabla f|(y) \le 1$ for $\meas$-a.e.\ $y \in X$ has $1$-Lipschitz representative.
\item{(Bochner inequality)} For any $f \in D(\Delta)$ with $\Delta f \in H^{1, 2}$,
\begin{gather*}
\frac{1}{2}\Delta |\nabla f|^2\ge \frac{(\Delta f)^2}{N}+\langle \nabla \Delta f, \nabla f\rangle +K|\nabla f|^2
\end{gather*}
holds in the weak sense, that is,
\begin{gather*}
\frac{1}{2}\int_X\Delta \phi |\nabla f|^2\di \meas \ge \int_X \phi \left(\frac{(\Delta f)^2}{N}+\langle \nabla \Delta f, \nabla f\rangle +K|\nabla f|^2\right)\di \meas
\end{gather*}
for any $\phi \in D(\Delta) \cap L^{\infty}(X, \meas)$ with $\Delta \phi \in L^{\infty}(X, \meas)$ and $\phi \ge 0$, where
\begin{align*}
D(\Delta) :=\bigg\{& f \in H^{1, 2}; \, \exists\, h=:\Delta f \in L^2, \\
&\mathrm{s.t.} \, \int_X\langle \nabla f, \nabla \psi \rangle \di \meas
=-\int_Xh\psi \di \meas,\,\forall\, \psi \in H^{1, 2}\bigg\}.
\end{align*}
\end{itemize}
See also \cite{A, AmbrosioGigliMondinoRajala, AmbrosioMondinoSavare, CM}. It is worth pointing out that any Ricci limit space as obtained by (\ref{1}) is an $\RCD(K(n-1), n)$ space by the stability of $\RCD(K, N)$ spaces with respect to the pmGH convergence proved in \cite{GMS}, and that any $n$-dimensional Alexandrov space of curvature bounded below by $K$ with $\mathcal{H}^n$ is also an $\RCD(K(n-1), n)$ space, which is proved in \cite{Petrunin, ZZ}.

From now on we fix $K \in \mathbb{R}$ and a finite $N<\infty$.
Thanks to recent quick developments on the study of $\RCD(K, N)$ spaces, most of the well-known properties on Ricci limit spaces can be covered by the $\RCD$ theory. For example, it is proved in \cite{BrueSemola} that the essential dimension, denoted by $\dim(X, \dist, \meas)$,\footnote{We know that $\dim_{\mathcal{H}}(X, \dist) \ge \dim(X, \dist, \meas)$ holds. See the final section. It is also conjectured that $\dim_{\mathcal{H}}(X, \dist) = \dim(X, \dist, \meas)$ holds. See also~\cite{Naber}.} whose definition is the unique $k$ such that the $k$-dimensional regular set has positive $\meas$-measure, is well-defined. This gives a generalization of a result proved in \cite{CN} to $\RCD(K, N)$ spaces. See also~\cite{Kit}.

On the other hand, a special class of $\RCD(K, N)$ spaces, so-called \textit{non-collapsed $\RCD(K, N)$ spaces}, is proposed in \cite{DePhillippisGigli} as the synthetic counterpart of non-collapsed Ricci limit spaces. An $\RCD (K, N)$ space is said to be \textit{non-collapsed} if $\meas=\mathcal{H}^N$ holds. Let us emphasize that it is essential that the upper bound $N$ of dimension (as $\RCD$ spaces) coincides with the dimension~$N$ of the Hausdorff measure $\mathcal{H}^N$ in the definition of the non-collapsed condition.\footnote{In general, the optimal dimension $\dim_{\mathrm{opt}}(X, \dist, \meas)$ of $(X, \dist, \meas)$ as $\RCD$ spaces, defined by the infimum of $\hat{N}$ such that $(X, \dist, \meas)$ is an $\RCD\big(\hat{K}, \hat{N}\big)$ space for some $\hat{K}$, is not equal to $\dim_{\mathcal{H}}(X, \dist)$. It follows from a result proved in \cite{DePhillippisGigli} that $[\dim_{\mathrm{opt}}(X, \dist, \meas)] \ge \dim_{\mathcal{H}}(X, \dist)$ holds, where~$[a]$ denotes the integer part of $a \in \mathbb{R}$. See also~\cite{HSZ}.} Then non-collapsed $\RCD(K, N)$ spaces have nicer properties, including that for non-collapsed Ricci limit spaces, rather than general $\RCD(K, N)$ spaces. For instance, any $N$-dimensional regular point~$p$ has a neighbourhood $U_p$ of~$p$ which is homeomorphic to~$\mathbb{R}^N$ (see also~\cite{KM}).

Let us mention that $\dim(X, \dist, \meas)$ is equal to $N$ if $(X, \dist, \meas)$ is a non-collapsed $\RCD(K, N)$ space. It is conjectured in~\cite{DePhillippisGigli} that the converse implication is also true up to multiplication by a positive constant to the measure, that is:
\begin{Conjecture}[De Philippis--Gigli]\label{16}
If an $\RCD(K, N)$ space $(X, \dist, \meas)$ satisfies $\dim(X, \dist, \meas)=N$, then $\meas=a\mathcal{H}^N$ holds for some $a \in (0, \infty)$.
\end{Conjecture}
This conjecture is true if $(X, \dist)$ is compact, which is proved in \cite{Honda}. Note that since the $\RCD(K, N)$ condition is unchanged under multiplication by a positive constant to the measure, if $\big(X, \dist, a\mathcal{H}^N\big)$ is an $\RCD(K, N)$ space, then $\big(X, \dist, \mathcal{H}^N\big)$ is a non-collapsed $\RCD(K, N)$ space.

Conjecture \ref{4} can be formulated in the $\RCD$ setting as follows:
\begin{Conjecture}\label{6}
If an $\RCD(K, N)$ space $(X, \dist, \meas)$ satisfies
\begin{gather}\label{7}
\meas=b\mathcal{H}^k
\end{gather} for some $b, k \in (0, \infty)$, then $\big(X, \dist, \mathcal{H}^k\big)$ is a non-collapsed $\RCD(K, k)$ space.
\end{Conjecture}
Note that the implication from Conjectures \ref{6} to \ref{4} is trivial by letting $n=N$.
Combining a result proved in \cite{BrueSemola} with Conjectures \ref{16} and \ref{6}, we also propose:
\begin{Conjecture}\label{102}
Let $(X, \dist, \meas)$ be an $\RCD (K, N)$ space. Then the following two conditions $(a)$ and $(b)$ are equivalent:
\begin{enumerate}\itemsep=0pt
\item[$(a)$] For all $f \in D(\Delta)$,
\begin{gather}\label{100}
\Delta f(x)=\operatorname{tr}(\mathrm{Hess}_f)(x)
\end{gather}
holds for $\meas$-a.e.\ $x \in X$, where the Hessian, $\mathrm{Hess}_f$, of $f$ is defined in~{\rm \cite{Gigli}} as a $(0, 2)$-type~$L^2$ tensor.
\item[$(b)$] \eqref{7} holds for some $b, k \in (0, \infty)$.
\end{enumerate}
\end{Conjecture}
One implication, from (a) to (b), is equivalent to Conjecture \ref{16}. Let us check this fact.

Assume that Conjecture \ref{16} is true. If (a) holds, then
it follows from a result of \cite{BrueSemola}, which confirms a conjecture raised in \cite{DePhillippisGigli}, that $(X, \dist, \meas)$ is an $\RCD(K, k)$ space, where $k=\dim(X, \dist, \meas)$. In particular, Conjecture \ref{16} yields (b).
Thus we have the implication from (a) to (b).

Next assume that the implication from (a) to (b) holds. If an $\RCD(K, N)$ space $(X, \dist, \meas)$ satisfies $\dim(X, \dist, \meas)=N$, then a result proved in \cite{Han} shows that (a) holds. Thus we have (b), which implies that Conjecture \ref{16} is true because (\ref{7}) implies $k=\dim(X, \dist, \meas)$ (see the final section).

On the other hand, the other implication, from (b) to (a), follows from Conjecture \ref{6}. That is, if (b) holds, then Conjecture \ref{6} yields that $\big(X, \dist, \mathcal{H}^k\big)$ is a non-collapsed $\RCD(K, k)$ space. In particular applying a result of \cite{Han} again implies (a).

It is worth pointing out that for a compact $\RCD(K, N)$ space $\big(X, \dist, \mathcal{H}^k\big)$, it is proved in~\cite{Honda} that $\big(X, \dist, \mathcal{H}^k\big)$ is a non-collapsed $\RCD(K, k)$ space if and only if
\begin{gather}\label{17}
\inf_{x \in X, r \in (0, 1)}\frac{\mathcal{H}^k(B_r(x))}{r^k}>0
\end{gather}
holds. Thus in the case when $(X, \dist)$ is compact, the remaining issue for Conjecture \ref{102} is only to prove (\ref{17}) under assuming (\ref{7}).

Finally let us go back to the question $(\heartsuit)$ appeared in the previous section. It is constructed in \cite{CheegerColding1} that for any sufficiently small $\epsilon>0$ there exists a sequence of complete Riemannian metrics $g_i$ on $\mathbb{R}^8$ with positive Ricci curvature such that
\begin{gather*}
\left(\mathbb{R}^8, \dist_{g_i}, 0, \frac{\mathrm{vol}_{g_i}}{\mathrm{vol}_{g_i}B_1(0)} \right) \stackrel{\mathrm{pmGH}}{\longrightarrow} \big(\mathbb{R}_{\ge 0} \times \mathbb{S}^4, \dist_{h_{\epsilon}}, p, \nu_{\epsilon}\big), \qquad h_{\epsilon}=(\di t)^2 +\left(\frac{t^{1+\epsilon}}{2}\right)^2g_{\mathbb{S}^4}
\end{gather*}
holds for some Borel measure $\nu_{\epsilon}$, where $p$ is the cusp and $g_{\mathbb{S}^4}$ is the standard Riemannian metric on $\mathbb{S}^4$. This limit space is called a metric horn. Let us denote it by $(X_{\epsilon}, \dist_{\epsilon}, p, \meas_{\epsilon})$.

From now on we will check that for any $K_1 \in \mathbb{R}$, $\big(X_{\epsilon}, \dist_{\epsilon}, \mathcal{H}^5\big)$ is not an $\RCD(K_1, \infty)$ space (in particular, $\big(X_{\epsilon}, \dist_{\epsilon}, \mathcal{H}^5\big)$ cannot be an $\RCD(K_1, N_1)$ space for any $N_1 \ge 5$). Before proving it, let us remark by definition of $h_{\epsilon}$ that we have
\begin{gather}\label{800}
\lim_{r \to 0^+}\frac{\mathcal{H}^5(B_r(p))}{r^5}=0.
\end{gather}

Assume that $\big(X_{\epsilon}, \dist_{\epsilon}, \mathcal{H}^5\big)$ is an $\RCD(K_1, \infty)$ space for some $K_1 \in \mathbb{R}$. Note that it is easy to see that an open set $B_R(p) \setminus \overline{B}_r(p)$ for all $r, R \in (0, \infty)$ with $r<R$ is isometric as Riemannian manifolds to an open subset of a closed Riemannian manifold (by gluing ``two caps'' along the boundary). In particular the localities of the minimal relaxed slope, of the Laplacian and of the Hessian (see~\cite{Gigli}) yield
\begin{gather*}
\int_{X_{\epsilon}}\langle \nabla f, \nabla \phi \rangle\di \mathcal{H}^5=-\int_{X_{\epsilon}}\operatorname{tr}(\mathrm{Hess}_f) \phi \di \mathcal{H}^5, \qquad \forall\, f \in D(\Delta),\quad \forall\, \phi \in C^{\infty}_c(X_{\epsilon} \setminus \{p\}).
\end{gather*}
Thus we see that $\Delta f(x) = \operatorname{tr}(\mathrm{Hess}_f)(x)$ for $\mathcal{H}^5$-a.e.\ $x \in X_{\epsilon}$. Then it follows from the Bochner inequality involving the Hessian term proved in \cite{Gigli} that in the weak sense:
\begin{align*}
\frac{1}{2}\Delta |\nabla f|^2 &\ge |\mathrm{Hess}_{f}|^2 +\langle \nabla \Delta f, \nabla f\rangle +K_1|\nabla f|^2 \nonumber \\
&\ge \frac{(\operatorname{tr}(\mathrm{Hess}_{f}))^2}{5}+\langle \nabla \Delta f, \nabla f\rangle +K_1|\nabla f|^2 \nonumber \\
&=\frac{(\Delta f)^2}{5}+\langle \nabla \Delta f, \nabla f\rangle +K_1|\nabla f|^2
\end{align*}
holds for any $f \in D(\Delta)$ with $\Delta f \in H^{1, 2}$.
This shows that $\big(X_{\epsilon}, \dist_{\epsilon}, \mathcal{H}^5\big)$ is an $\RCD(K_1, 5)$ space, that is, it is non-collapsed. In particular the Bishop--Gromov inequality yields
\begin{gather*}
\liminf_{r \to 0^+}\frac{\mathcal{H}^5(B_r(p))}{r^5}>0,
\end{gather*}
which contradicts (\ref{800}).

\section{Regularity on reference measure}
Let $(X, \dist, \meas)$ be an $\RCD(K, N)$ space with $k=\dim(X, \dist, \meas)$.
We recall some structure results on the reference measure $\meas$.
Based on results obtained in \cite{ DePhillippisMarcheseRindler, GigliPasqualetto, KellMondino, MondinoNaber}, it is proved in \cite{AmbrosioHondaTewodrose} that the limit
\begin{gather}\label{9}
\phi(x):=\lim_{r \to 0^+}\frac{\meas (B_r(x))}{\omega_kr^k}
\end{gather}
exists in $(0, \infty)$ for $\meas$-a.e.\ $x \in X$ (we denote by $\mathcal{R}_k^*$ the set of all such points $x\in X$) and coincides with the Radon-Nikodym derivative of the restriction of $\meas$ to $\mathcal{R}_k^*$ with respect to $\mathcal{H}^k$, where $\omega_k=\mathcal{H}^k(B_1(0_k))$.
In particular we have $\dim_{\mathcal{H}}(X, \dist) \ge k$ because $\mathcal{H}^k(\mathcal{R}_k^*)>0$ holds, which comes from the finiteness of $\phi$ on $\mathcal{R}_k^*$ with $\meas (\mathcal{R}_k^*)>0$.
Moreover if $\meas=b\mathcal{H}^l$ holds for some $b, l \in \mathbb{R}_{>0}$, then $l=k$ holds because $\mathcal{H}^l$ and $\mathcal{H}^k$ are mutually absolutely continuous on $\mathcal{R}_k^*$.

Let $f$ be a function defined on a Borel subset $A$ of $X$.
We say that $f$ is differentiable for $\meas$-a.e.\ $x \in A$ if there exists a family of Borel subset $\{A_i\}_{i \in \mathbb{N}}$ of $A$ such that $\meas (A \setminus \bigcup_iA_i)=0$ holds and that $f|_{A_i}$ is Lipschitz for any $i$. See \cite{ABT, Honda1}. Typical examples can be found in Sobolev functions, that is, we see that any $g \in H^{1, p}(U, \dist, \meas)$ for some open subset $U$ of $X$ and some $p \in [1, \infty)$ is differentiable for $\meas$-a.e.\ $x \in U$ via the standard telescopic argument.\footnote{The characteristic function $\chi_{B_1(0_n)}$ of $B_1(0_n)$ in $\mathbb{R}^n$ is differentiable for $\mathcal{H}^n$-a.e.\ $x \in \mathbb{R}^n$, but for all $p \in [1, \infty)$ it is not in $H^{1, p}$.} Let us prove this fact for reader's convenience.

Fix $x \in U$ and find $R>0$ with $B_{3R}(x) \subset U$. For any $L>1$ let us define:
\begin{gather*}
A_L:=\big\{y \in B_R(x);\, \big| |\nabla g|_p^p \big|_{B_r(y)} \le L^p, \,\forall\, r \in (0, R]\big\},
\end{gather*}
where $f_A$ is the average of $f$ over $A \subset X$ and $|\nabla g|_p$ denotes the minimal ($p$-)weak upper gradient of $g$ (see~\cite{GHan} for the independence on~$p$ of $|\nabla g|_p$). Note that $\meas (B_R(x) \setminus A_L) \to 0$ holds as $L \uparrow \infty$ by the maximal function theorem. Thus it is enough to check:
\begin{enumerate}\itemsep=0pt
\item[$(\clubsuit)$] there exists a Borel subset $\hat{A}_L$ of $A_L$ such that $\meas \big(A_L \setminus \hat{A}_L\big)=0$ holds and that $g|_{\hat{A}_L}$ is Lipschitz.
\end{enumerate}

The $(1, p)$-Poincar\'e inequality proved in \cite{Raj} (for more general class, $\mathrm{CD}(K, \infty)$ spaces), with the volume doubling property which follows from the Bishop--Gromov inequality yields
\begin{gather*}
\big| g- g_{B_r(y)}\big|_{B_r(y)} \le C(K, N)Lr, \qquad \forall\, y \in A_L,\quad \forall\, r \in (0, R].
\end{gather*}
See also \cite{HK}. In particular for any $i \in \mathbb{Z}_{\ge 0}$, any $r \in [0, R]$ and any $y \in A_L$, we have
\begin{align*}
\big| g_{B_{2^{-(i+1)}r}(y)}-g_{B_{2^{-i}r}(y)}\big| &\le \big| g- g_{B_{2^{-i}r}(y)}\big|_{B_{2^{-(i+1)}r}(y)} \nonumber \\
&\le C(K, N, R) \big| g- g_{B_{2^{-i}r}(y)}\big|_{B_{2^{-i}r}(y)} \le C(K, N, R)2^{-i}Lr.
\end{align*}
Thus taking the sum with respect to $i$ yields that any Lebesgue point $y \in A_L$ of $g$ satisfies
\begin{gather}\label{eqeq}
 \big|g(y)-g_{B_r(y)} \big| \le C(K, N, R)Lr.
\end{gather}
Since it is easy to see $|g_{B_r(y)}-g_{B_r(z)}|\le C(K, N, R)L$ for all $y, z \in A_L$ with $\dist (y, z) \le r$, (\ref{eqeq})~shows
\begin{gather}\label{eee}
|g(y)-g(z)|\le C(K, N, R)Lr
\end{gather}
for all Lebesgue points $y, z \in A_L$ of $g$ with $\dist (y, z) \le r$.
In particular taking $r=\dist(y, z)$ in~(\ref{eee}) completes the proof of $(\clubsuit)$.

It is trivial from \cite{Gigli} that if $f$ is differentiable for $\meas$-a.e.\ $x \in A$ with $\meas (X \setminus A)=0$, then $\nabla f$ is well-defined in $L^0(TX)$, that is, $\nabla f$ is a Borel measurable vector field on~$X$.

We are now in a position to introduce the final conjecture (recall $\meas ( X \setminus \mathcal{R}_k^*)=0$, where $k=\dim(X, \dist, \meas)$):
\begin{Conjecture}\label{10} The function $\phi$ defined by \eqref{9} is differentiable for $\meas$-a.e.\ $x \in \mathcal{R}_k^*$. Moreover for all $f \in D(\Delta)$,
\begin{gather*}
\Delta f(x)=\operatorname{tr}(\mathrm{Hess}_f)(x)+\langle \nabla \log \phi, \nabla f\rangle (x)
\end{gather*}
holds for $\meas$-a.e.\ $x \in X$.
\end{Conjecture}
Let us remark that the implication from Conjecture \ref{10} to (\ref{100}) is trivial. A partial contribution to Conjecture \ref{10} can be found in \cite{Honda}.

The technique provided in \cite{Honda} is useful for all conjectures above in the case when $(X, \dist)$ is compact.
This is to apply a geometric flow defined by embedding maps in $L^2$ via the global heat kernel.\footnote{The Dirichlet heat kernel on a bounded open subset also works along this direction. See also \cite{ZZ1}.} Such embedding maps are introduced and studied first in \cite{BerardBessonGallot} for closed Riemannian manifolds. Recently in \cite{AHPT} this observation is generalized to $\RCD(K, N)$ spaces by using stability results of Sobolev functions with repect to the pmGH convergence proved in \cite{AmbrosioHonda, AmbrosioHonda2}. It seems to the author that this technique, using a geometric flow, is useful for all conjectures proposed in this paper even in the case when $(X, \dist)$ is non-compact.

\subsection*{Acknowledgements}

The author acknowledges supports of the Grantin-Aid for Young Scientists (B) 16K17585 and Grant-in-Aid for Scientific Research (B) of 18H01118. He would like to express his appreciation to the referees for valuable suggestions on the first version, which make the paper more readable.


\pdfbookmark[1]{References}{ref}
\LastPageEnding


\begin{thebibliography}{99}
\footnotesize\itemsep=0pt

\bibitem{A}
Ambrosio L., Calculus, heat flow and curvature-dimension bounds in metric
 measure spaces, in Proceedings of the {I}nternational {C}ongress of
 {M}athematicians~-- {R}io de {J}aneiro 2018, {V}ol.~{I}, {P}lenary lectures,
 World Sci. Publ., Hackensack, NJ, 2018, 301--340.

\bibitem{ABT}
Ambrosio L., Bru\`e E., Trevisan D., Lusin-type approximation of {S}obolev by
 {L}ipschitz functions, in {G}aussian and {${\rm RCD}(K,\infty)$} spaces,
 \href{https://doi.org/10.1016/j.aim.2018.09.033}{\textit{Adv. Math.}} \textbf{339} (2018), 426--452, \href{https://arxiv.org/abs/1712.06315}{arXiv:1712.06315}.

\bibitem{AmbrosioGigliMondinoRajala}
Ambrosio L., Gigli N., Mondino A., Rajala T., Riemannian {R}icci curvature
 lower bounds in metric measure spaces with {$\sigma$}-finite measure,
 \href{https://doi.org/10.1090/S0002-9947-2015-06111-X}{\textit{Trans. Amer. Math. Soc.}} \textbf{367} (2015), 4661--4701, \href{https://arxiv.org/abs/1207.4924}{arXiv:1207.4924}.

\bibitem{AmbrosioGigliSavare14}
Ambrosio L., Gigli N., Savar\'{e} G., Metric measure spaces with {R}iemannian
 {R}icci curvature bounded from below, \href{https://doi.org/10.1215/00127094-2681605}{\textit{Duke Math.~J.}} \textbf{163} (2014), 1405--1490, \href{https://arxiv.org/abs/1109.0222}{arXiv:1109.0222}.

\bibitem{AmbrosioHonda}
Ambrosio L., Honda S., New stability results for sequences of metric measure
 spaces with uniform {R}icci bounds from below, in Measure Theory in
 Non-Smooth Spaces, Partial Differ. Equ. Meas. Theory, \href{https://doi.org/10.1515/9783110550832-001}{De Gruyter Open}, Warsaw, 2017, 1--51, \href{https://arxiv.org/abs/1605.07908}{arXiv:1605.07908}.

\bibitem{AmbrosioHonda2}
Ambrosio L., Honda S., Local spectral convergence in {${\rm RCD}^*(K,N)$}
 spaces, \href{https://doi.org/10.1016/j.na.2017.04.003}{\textit{Nonlinear Anal.}} \textbf{177} (2018), 1--23, \href{https://arxiv.org/abs/1703.04939}{arXiv:1703.04939}.

\bibitem{AHPT}
Ambrosio L., Honda S., Portegies J.W., Tewodrose D., Embedding of ${\rm
 RCD}^*(K, N)$-spaces in $L^2$ via eigenfunctions, \href{https://arxiv.org/abs/1812.03712}{arXiv:1812.03712}.

\bibitem{AmbrosioHondaTewodrose}
Ambrosio L., Honda S., Tewodrose D., Short-time behavior of the heat kernel and
 {W}eyl's law on {${\rm RCD}^*(K,N)$} spaces, \href{https://doi.org/10.1007/s10455-017-9569-x}{\textit{Ann. Global Anal. Geom.}}
 \textbf{53} (2018), 97--119, \href{https://arxiv.org/abs/1701.03906}{arXiv:1701.03906}.

\bibitem{AmbrosioMondinoSavare}
Ambrosio L., Mondino A., Savar\'{e} G., Nonlinear diffusion equations and
 curvature conditions in metric measure spaces, \href{https://doi.org/10.1090/memo/1270}{\textit{Mem. Amer. Math. Soc.}} \textbf{262} (2019), v+121~pages, \href{https://arxiv.org/abs/1509.07273}{arXiv:1509.07273}.

\bibitem{An}
Anderson M.T., Hausdorff perturbations of {R}icci-flat manifolds and the
 splitting theorem, \href{https://doi.org/10.1215/S0012-7094-92-06803-7}{\textit{Duke Math.~J.}} \textbf{68} (1992), 67--82.

\bibitem{BerardBessonGallot}
B\'{e}rard P., Besson G., Gallot S., Embedding {R}iemannian manifolds by their
 heat kernel, \href{https://doi.org/10.1007/BF01896401}{\textit{Geom. Funct. Anal.}} \textbf{4} (1994), 373--398.

\bibitem{BrueSemola}
Bru\`e E., Semola D., Constancy of dimension for ${\rm RCD}^*(K,N)$ spaces via
 regularity of {L}agrangian flows, \href{https://doi.org/10.1002/cpa.21849}{\textit{Comm. Pure Appl. Math.}}, {t}o appear, \href{https://arxiv.org/abs/1803.04387}{arXiv:1803.04387}.

\bibitem{BGP}
Burago Y., Gromov M., Perelman G., A.{D}.~{A}leksandrov spaces with curvatures
 bounded below, \href{https://doi.org/10.1070/RM1992v047n02ABEH000877}{\textit{Russian Math. Surveys}} \textbf{47} (1992), no.~2, 1--58.

\bibitem{CM}
Cavalletti F., Milman E., The globalization theorem for the curvature dimension
 condition, \href{https://arxiv.org/abs/1612.07623}{arXiv:1612.07623}.

\bibitem{CheegerColding1}
Cheeger J., Colding T.H., On the structure of spaces with {R}icci curvature
 bounded below.~{I}, \href{https://doi.org/10.4310/jdg/1214459974}{\textit{J.~Differential Geom.}} \textbf{46} (1997), 406--480.

\bibitem{CheegerColding2}
Cheeger J., Colding T.H., On the structure of spaces with {R}icci curvature
 bounded below.~{II}, \href{https://doi.org/10.4310/jdg/1214342145}{\textit{J.~Differential Geom.}} \textbf{54} (2000), 13--35.

\bibitem{CheegerColding3}
Cheeger J., Colding T.H., On the structure of spaces with {R}icci curvature
 bounded below.~{III}, \href{https://doi.org/10.4310/jdg/1214342146}{\textit{J.~Differential Geom.}} \textbf{54} (2000), 37--74.

\bibitem{CJN}
Cheeger J., Jiang W., Naber A., Rectifiability of singular sets in noncollapsed
 spaces with {R}icci curvature bounded below, \href{https://arxiv.org/abs/1805.07988}{arXiv:1805.07988}.

\bibitem{CN}
Colding T.H., Naber A., Sharp {H}\"{o}lder continuity of tangent cones for
 spaces with a lower {R}icci curvature bound and applications, \href{https://doi.org/10.4007/annals.2012.176.2.10}{\textit{Ann. of
 Math.}} \textbf{176} (2012), 1173--1229, \href{https://arxiv.org/abs/1102.5003}{arXiv:1102.5003}.

\bibitem{CN1}
Colding T.H., Naber A., Characterization of tangent cones of noncollapsed
 limits with lower {R}icci bounds and applications, \href{https://doi.org/10.1007/s00039-012-0202-7}{\textit{Geom. Funct.
 Anal.}} \textbf{23} (2013), 134--148, \href{https://arxiv.org/abs/1108.3244}{arXiv:1108.3244}.

\bibitem{DePhillippisGigli}
De~Philippis G., Gigli N., Non-collapsed spaces with {R}icci curvature bounded
 from below, \href{https://doi.org/10.5802/jep.80}{\textit{J.~\'Ec. polytech. Math.}} \textbf{5} (2018), 613--650,
 \href{https://arxiv.org/abs/1708.02060}{arXiv:1708.02060}.

\bibitem{DePhillippisMarcheseRindler}
De~Philippis G., Marchese A., Rindler F., On a conjecture of {C}heeger, in
 Measure Theory in Non-Smooth Spaces, \textit{Partial Differ. Equ. Meas. Theory}, \href{https://doi.org/10.1515/9783110550832-004}{De
 Gruyter Open}, Warsaw, 2017, 145--155, \href{https://arxiv.org/abs/1607.02554}{arXiv:1607.02554}.

\bibitem{ErbarKuwadaSturm}
Erbar M., Kuwada K., Sturm K.-T., On the equivalence of the entropic
 curvature-dimension condition and {B}ochner's inequality on metric measure
 spaces, \href{https://doi.org/10.1007/s00222-014-0563-7}{\textit{Invent. Math.}} \textbf{201} (2015), 993--1071,
 \href{https://arxiv.org/abs/1303.4382}{arXiv:1303.4382}.

\bibitem{F}
Fukaya K., Collapsing of {R}iemannian manifolds and eigenvalues of {L}aplace
 operator, \href{https://doi.org/10.1007/BF01389241}{\textit{Invent. Math.}} \textbf{87} (1987), 517--547.

\bibitem{Gigli13}
Gigli N., The splitting theorem in non-smooth context, \href{https://arxiv.org/abs/1302.5555}{arXiv:1302.5555}.

\bibitem{Gigli}
Gigli N., Nonsmooth differential geometry~-- an approach tailored for spaces
 with {R}icci curvature bounded from below, \href{https://doi.org/10.1090/memo/1196}{\textit{Mem. Amer. Math. Soc.}}
 \textbf{251} (2018), v+161~pages, \href{https://arxiv.org/abs/1407.0809}{arXiv:1407.0809}.

\bibitem{GHan}
Gigli N., Han B.-X., Independence on {$p$} of weak upper gradients on {${\rm
 RCD}$} spaces, \href{https://doi.org/10.1016/j.jfa.2016.04.014}{\textit{J.~Funct. Anal.}} \textbf{271} (2016), 1--11, \href{https://arxiv.org/abs/1407.7350}{arXiv:1407.7350}.

\bibitem{GMS}
Gigli N., Mondino A., Savar\'{e} G., Convergence of pointed non-compact metric
 measure spaces and stability of {R}icci curvature bounds and heat flows,
 \href{https://doi.org/10.1112/plms/pdv047}{\textit{Proc. Lond. Math. Soc.}} \textbf{111} (2015), 1071--1129, \href{https://arxiv.org/abs/1311.4907}{arXiv:1311.4907}.

\bibitem{GigliPasqualetto}
Gigli N., Pasqualetto E., Behaviour of the reference measure on ${\rm RCD}$
 spaces under charts, \textit{Comm. Pure Appl. Math.}, {t}o appear,
 \href{https://arxiv.org/abs/1607.05188}{arXiv:1607.05188}.

\bibitem{HK}
Haj{\l}asz P., Koskela P., Sobolev meets {P}oincar\'e, \textit{C.~R.~Acad. Sci.
 Paris S\'er.~I Math.} \textbf{320} (1995), 1211--1215.

\bibitem{Han}
Han B.-X., Ricci tensor on {${\rm RCD}^*(K,N)$} spaces, \href{https://doi.org/10.1007/s12220-017-9863-7}{\textit{J.~Geom. Anal.}}
 \textbf{28} (2018), 1295--1314, \href{https://arxiv.org/abs/1412.0441}{arXiv:1412.0441}.

\bibitem{H}
Hattori K., The nonuniqueness of the tangent cones at infinity of {R}icci-flat
 manifolds, \href{https://doi.org/10.2140/gt.2017.21.2683}{\textit{Geom. Topol.}} \textbf{21} (2017), 2683--2723, \href{https://arxiv.org/abs/1503.07278}{arXiv:1503.07278}.

\bibitem{Honda1}
Honda S., A weakly second-order differential structure on rectifiable metric
 measure spaces, \href{https://doi.org/10.2140/gt.2014.18.633}{\textit{Geom. Topol.}} \textbf{18} (2014), 633--668, \href{https://arxiv.org/abs/1112.0099}{arXiv:1112.0099}.

\bibitem{Honda}
Honda S., New differential operator and non-collapsed ${\rm RCD}$ spaces,
 \textit{Geom. Topol.}, {t}o appear, \href{https://arxiv.org/abs/1905.00123}{arXiv:1905.00123}.

\bibitem{HSZ}
Honda S., Sun S., Zhang R., A note on the collapsing geometry of
 hyperk\"{a}hler four manifolds, \href{https://doi.org/10.1007/s11425-019-1602-x}{\textit{Sci. China Math.}} \textbf{62} (2019),
 2195--2210.

\bibitem{K}
Kapovitch V., Perelman's stability theorem, in Surveys in Differential
 Geometry, {V}ol.~{XI}, \textit{Surv. Differ. Geom.}, Vol.~11, \href{https://doi.org/10.4310/SDG.2006.v11.n1.a5}{Int. Press},
 Somerville, MA, 2007, 103--136, \href{https://arxiv.org/abs/math.DG/0703002}{arXiv:math.DG/0703002}.

\bibitem{KM}
Kapovitch V., Mondino A., On the topology and the boundary of $N$-dimensional
 ${\rm RCD}(K,N)$ spaces, \href{https://arxiv.org/abs/11907.02614}{arXiv:11907.02614}.

\bibitem{KellMondino}
Kell M., Mondino A., On the volume measure of non-smooth spaces with {R}icci
 curvature bounded below, \textit{Ann. Sc. Norm. Super. Pisa Cl. Sci.~(5)}
 \textbf{18} (2018), 593--610, \href{https://arxiv.org/abs/1607.02036}{arXiv:1607.02036}.

\bibitem{Kit}
Kitabeppu Y., A sufficient condition to a regular set being of positive measure
 on {${\rm RCD}$} spaces, \href{https://doi.org/10.1007/s11118-018-9708-4}{\textit{Potential Anal.}} \textbf{51} (2019),
 179--196, \href{https://arxiv.org/abs/1708.04309}{arXiv:1708.04309}.

\bibitem{LottVillani}
Lott J., Villani C., Ricci curvature for metric-measure spaces via optimal
 transport, \href{https://doi.org/10.4007/annals.2009.169.903}{\textit{Ann. of Math.}} \textbf{169} (2009), 903--991,
 \href{https://arxiv.org/abs/math.DG/0412127}{arXiv:math.DG/0412127}.

\bibitem{M2}
Menguy X., Examples of nonpolar limit spaces, \href{https://doi.org/10.1353/ajm.2000.0041}{\textit{Amer.~J. Math.}} \textbf{122} (2000), 927--937.

\bibitem{M1}
Menguy X., Noncollapsing examples with positive {R}icci curvature and infinite
 topological type, \href{https://doi.org/10.1007/PL00001632}{\textit{Geom. Funct. Anal.}} \textbf{10} (2000), 600--627.

\bibitem{M3}
Menguy X., Examples of strictly weakly regular points, \href{https://doi.org/10.1007/PL00001667}{\textit{Geom. Funct.
 Anal.}} \textbf{11} (2001), 124--131.

\bibitem{MondinoNaber}
Mondino A., Naber A., Structure theory of metric measure spaces with lower
 {R}icci curvature bounds, \href{https://doi.org/10.4171/JEMS/874}{\textit{J.~Eur. Math. Soc.}} \textbf{21} (2019), 1809--1854, \href{https://arxiv.org/abs/1405.2222}{arXiv:1405.2222}.

\bibitem{Naber}
Naber A., The geometry of {R}icci curvature, in Proceedings of the
 {I}nternational {C}ongress of {M}athematicians~-- {S}eoul 2014, {V}ol.~{II},
 Kyung Moon Sa, Seoul, 2014, 911--937.

\bibitem{P}
Perelman G., Spaces with curvature bounded below~II, {P}reprint.

\bibitem{Petrunin}
Petrunin A., Alexandrov meets {L}ott--{V}illani--{S}turm,
 \textit{M\"{u}nster~J. Math.} \textbf{4} (2011), 53--64, \href{https://arxiv.org/abs/1003.5948}{arXiv:1003.5948}.

\bibitem{Raj}
Rajala T., Local {P}oincar\'{e} inequalities from stable curvature conditions
 on metric spaces, \href{https://doi.org/10.1007/s00526-011-0442-7}{\textit{Calc. Var. Partial Differential Equations}}
 \textbf{44} (2012), 477--494, \href{https://arxiv.org/abs/1107.4842}{arXiv:1107.4842}.

\bibitem{Sturm06a}
Sturm K.-T., On the geometry of metric measure spaces.~{I}, \href{https://doi.org/10.1007/s11511-006-0002-8}{\textit{Acta Math.}}
 \textbf{196} (2006), 65--131.

\bibitem{Sturm06b}
Sturm K.-T., On the geometry of metric measure spaces.~{II}, \href{https://doi.org/10.1007/s11511-006-0003-7}{\textit{Acta Math.}}
 \textbf{196} (2006), 133--177.

\bibitem{Villani}
Villani C., Optimal transport: old and new, \textit{Grundlehren der Mathematischen
 Wissenschaften}, Vol.~338, \href{https://doi.org/10.1007/978-3-540-71050-9}{Springer-Verlag}, Berlin, 2009.
\bibitem{Y}
Yamaguchi T., Collapsing and pinching under a lower curvature bound,
 \href{https://doi.org/10.2307/2944340}{\textit{Ann. of Math.}} \textbf{133} (1991), 317--357.

\bibitem{ZZ}
Zhang H.-C., Zhu X.-P., Ricci curvature on {A}lexandrov spaces and rigidity
 theorems, \href{https://doi.org/10.4310/CAG.2010.v18.n3.a4}{\textit{Comm. Anal. Geom.}} \textbf{18} (2010), 503--553,
 \href{https://arxiv.org/abs/0912.3190}{arXiv:0912.3190}.

\bibitem{ZZ1}
Zhang H.-C., Zhu X.-P., Weyl's law on ${\rm RCD}^*(K, N)$ metric measure spaces,
\href{https://dx.doi.org/10.4310/CAG.2019.v27.n8.a8}{\textit{Comm. Anal. Geom.}} \textbf{27} (2019), 1869--1914, \href{https://arxiv.org/abs/1701.01967}{arXiv:1701.01967}.

\end{thebibliography}
\end{document}